\documentclass[reqno]{amsart}

\usepackage{amssymb}

\usepackage{euler,eucal}

\DeclareMathAlphabet{\mathbf}{T1}{ppl}{bx}{n}
\DeclareMathAlphabet{\mathrm}{T1}{ppl}{m}{n}

\numberwithin{equation}{section}

\newcommand\note[1]%
{$^\dagger$\marginpar{\footnotesize{$^\dagger${#1}}}}

\def\({\left(}
\def\){\right)}
\newcommand\eps{\varepsilon}

\newtheorem{theorem}{Theorem}[section]
\newtheorem{proposition}[theorem]{Proposition}
\newtheorem{lemma}[theorem]{Lemma}

\newtheorem{corollary}[theorem]{Corollary}

\theoremstyle{definition}

\newtheorem{example}[theorem]{Example}

\newcommand\lie{\mathfrak}

\newcommand\g{\lie{g}}

\newcommand\bb[1]{{\text{\bf#1}}}

\newcommand\R{\bb{R}}

\newcommand\ca{\mathcal}

\newcommand\func[1]{\operatorname{\mathrm{#1}}}
\newcommand\funclim[1]{\operatorname*{\mathrm{#1}}}

\newcommand\Ad{\func{Ad}}

\renewcommand\det{\func{det}}

\renewcommand\exp{\func{exp}}

\newcommand\id{\func{id}}
\newcommand\im{\func{im}}

\renewcommand\ker{\func{ker}}
\renewcommand\lim{\funclim{lim}}

\newcommand\inj{\hookrightarrow}
\newcommand\sur{\mathrel{\to\kern-1.8ex\to}}
\newcommand\iso{\mathrel{\hookrightarrow\kern-1.8ex\to}}

\newcommand\longhookrightarrow{\lhook\joinrel\longrightarrow}

\newcommand\longsur{\mathrel{\longrightarrow\kern-1.8ex\to}}
\newcommand\longiso{\mathrel{\longhookrightarrow\kern-1.8ex\to}}

\renewcommand\subset{\subseteq}

\begin{document}

%%%%%%%%%%%%%%%%%%%%%%%%%%%%%%%%%%%%%%%%%%%%%%%%%%%%%%%%%%%%%%%%%%%%%%%%
%%%%%%%%%%%%%%%%%%%%%%%%%%%%%%%%%%%%%%%%%%%%%%%%%%%%%%%%%%%%%%%%%%%%%%%%

\title{Equivariant symplectic Hodge theory and the $d_G\delta$-lemma}

\author{Yi Lin}

\email{linyi@math.cornell.edu}

\author{Reyer Sjamaar}

\email{sjamaar@math.cornell.edu}

\address{Department of Mathematics, Cornell University, Ithaca, New
York 14853-7901, USA}

\thanks{R. Sjamaar was partially supported by NSF Grant DMS-0071625.}

\date{15 September 2003}

\subjclass[2000]{53D20, 58A12}

\begin{abstract}
Consider a Hamiltonian action of a compact Lie group on a symplectic
manifold which has the strong Lefschetz property.  We establish
an equivariant version of the Merkulov-Guillemin $d\delta$-lemma and
an improved version of the Kirwan-Ginzburg equivariant formality
theorem, which says that every cohomology class has a \emph{canonical}
equivariant extension.
\end{abstract}

\maketitle

%%%%%%%%%%%%%%%%%%%%%%%%%%%%%%%%%%%%%%%%%%%%%%%%%%%%%%%%%%%%%%%%%%%%%%%%
%%%%%%%%%%%%%%%%%%%%%%%%%%%%%%%%%%%%%%%%%%%%%%%%%%%%%%%%%%%%%%%%%%%%%%%%

%%%%%%%%%%%%%%%%%%%%%%%%%%%%%%%%%%%%%%%%%%%%%%%%%%%%%%%%%%%%%%%%%%%%%%%%
\section{Introduction}
%%%%%%%%%%%%%%%%%%%%%%%%%%%%%%%%%%%%%%%%%%%%%%%%%%%%%%%%%%%%%%%%%%%%%%%%

Mimicking Riemannian Hodge theory Brylinski
\cite{brylinski;differential-poisson} introduced the symplectic Hodge
star operator and defined a differential form $\alpha$ on a symplectic
manifold to be \emph{harmonic} if $d\alpha=d{*}\alpha=0$.  He
conjectured that on a compact symplectic manifold every de Rham
cohomology class contains a harmonic representative and proved this
conjecture to be true for K\"ahler manifolds and certain other
examples.

A compact symplectic manifold $(M,\omega)$ of dimension $2n$ is said
to have the \emph{strong Lefschetz property} if the map
$$H^{n-k}(M)\longrightarrow H^{n+k}(M),\qquad c\longmapsto
[\omega]^k\wedge c$$
is an isomorphism for each $0\leq k \leq n$.  Mathieu
\cite{mathieu;harmonic-symplectic} proved the remarkable theorem that
Brylinski's conjecture is true for $M$ if and only if it has the
strong Lefschetz property.  This result was sharpened by Merkulov
\cite{merkulov;formality-canonical-symplectic} and Guillemin
\cite{guillemin;symplectic-hodge}, who independently established the
symplectic $d\delta$-lemma.  Let $\delta=\pm{*}d{*}$ be Koszul's
boundary operator and suppose $M$ has the strong Lefschetz property.
The $d\delta$-lemma asserts that if $\alpha$ is a harmonic $k$-form on
$M$ that is either exact or coexact then $\alpha=d\delta\beta$ for
some $k$-form $\beta$.

The object of this paper is to extend these results to an equivariant
setting.  Let $M$ be a compact symplectic manifold satisfying the
strong Lefschetz condition and let $G$ be a compact connected Lie
group acting on $M$ in a Hamiltonian fashion.  One of our main results
is an equivariant version of the $d\delta$-lemma, namely the
$d_G\delta$-lemma, which is obtained by replacing $d$ with the
equivariant exterior derivative $d_G$ on the Cartan double complex
$\Omega_G(M)=(S\g^*\otimes\Omega(M))^G$.

It was proved by Kirwan \cite{kirwan;cohomology-quotients-symplectic}
and Ginzburg \cite{ginzburg;equivariant-cohomology-kahler} that the
spectral sequence of the Cartan complex associated to the horizontal
filtration degenerates at the first term.  This implies the
equivariant formality theorem, which states that the equivariant
cohomology group $H^*_G(M)$ is isomorphic to $(S\g^*)^G\otimes H^*(M)$
as a vector space.  This isomorphism depends on a choice of an
equivariantly closed extension for each closed form on $M$.  The
$d\delta$-lemma produces a very simple proof of the Kirwan-Ginzburg
theorem, which replaces the Morse theory in
\cite{kirwan;cohomology-quotients-symplectic,%
ginzburg;equivariant-cohomology-kahler} with Hodge theory and which
shows in addition that there is a \emph{canonical} choice for these
equivariantly closed extensions up to coboundaries and therefore a
\emph{canonical} isomorphism from $H^*_G(M)$ to $(S\g^*)^G\otimes
H^*(M)$.  The idea behind this result, which appears to be new even
for K\"ahler manifolds, goes back to Deligne
\cite{deligne;lefschetz-degenerescence}, who showed how Lefschetz type
properties lead to the collapse of spectral sequences.

Section \ref{section;symplectic-hodge} of this paper is a brief review
of symplectic Hodge theory, including the $d\delta$-lemma.  In Section
\ref{section;dgdelta} we deduce from the $d\delta$-lemma the strong
version of the equivariant formality theorem and the
$d_G\delta$-lemma.

%%%%%%%%%%%%%%%%%%%%%%%%%%%%%%%%%%%%%%%%%%%%%%%%%%%%%%%%%%%%%%%%%%%%%%%%
\section{Symplectic Hodge theory}\label{section;symplectic-hodge}
%%%%%%%%%%%%%%%%%%%%%%%%%%%%%%%%%%%%%%%%%%%%%%%%%%%%%%%%%%%%%%%%%%%%%%%%

Let $(V,\omega)$ be a $2n$-dimensional symplectic vector space.  The
symplectic form induces an isomorphism $\flat\colon V\to V^*$ defined
by $\langle a^\flat,b\rangle=\omega(a,b)$, where
$\langle\cdot,\cdot\rangle$ denotes the dual pairing between $V$ and
$V^*$.  Let $\sharp\colon V^*\to V$ be the inverse of $\flat$.  Then
$\omega(v,w)=\omega(v^\sharp,w^\sharp)$ is a symplectic form on $V^*$,
which induces bilinear pairings on $\Lambda V^*$ also denoted by
$\omega$ and defined by
$$
\omega(v_1\wedge\cdots\wedge v_k,w_1\wedge\cdots\wedge
w_k)=\det\omega(v_i,w_j)
$$
for $v_i$, $w_j\in V^*$.  It is easy to show that $\omega$ is
nondegenerate; furthermore on $\Lambda^kV^*$ it is anti-symmetric when
$k$ is odd and symmetric when $k$ is even.  For $a\in V$ let
$\iota(a)\colon\Lambda^kV^*\to\Lambda^{k-1}V^*$ denote inner product
by $a$ and for $u\in V^*$ let
$\eps(u)\colon\Lambda^kV^*\to\Lambda^{k+1}V^*$ denote exterior product
by $u$.  An easy calculation confirms the following fact.

\begin{lemma}\label{lemma;adjoint}
For all $u\in V^*$ the maps $\eps(u)$ and $-\iota(u^\sharp)$ are
adjoint with respect to $\omega$.
\end{lemma}

The symplectic Hodge star operator
$*\colon\Lambda^kV^*\to\Lambda^{2n-k}V^*$ is defined by
$$
u\wedge*v=\omega(u,v)\frac{\omega^n}{n!}
$$
for $u$, $v\in\Lambda^kV^*$.

\begin{example}
A direct calculation shows that on a $2n$-dimensional symplectic
vector space $*\omega/k!=\omega^{n-k}/(n-k)!$ for all $k$.  In short,
$*\exp\omega=\exp\omega$.
\end{example}

The Hodge star operator has the following properties.

\begin{lemma}[{\cite{brylinski;differential-poisson}}]
\label{lemma1.1}
\begin{enumerate}
\item\label{item;sum}
Let $(V_1,\omega_1)$ and $(V_2,\omega_2)$ be symplectic vector spaces
and let $(V,\omega)$ be their direct sum.  Denote the Hodge star
operators on $V_1$, $V_2$ and $V$ by $*_1$, $*_2$ and $*$,
respectively.  Then
$$*(u_1\wedge
u_2)=(-1)^{k_1k_2}(*_1u_1)\wedge(*_2u_2)=(*_2u_2)\wedge(*_1u_1)$$
for all $u_1\in\Lambda^{k_1}V_1^*$ and
$u_2\in\Lambda^{k_2}V_2^*$.
\item\label{item;**}
Let $(V,\omega)$ be a symplectic vector space.  Then $*(*u)=u$ for all
$u\in\Lambda^kV^*$, i.e.\ $*^2=\id$.
\end{enumerate}
\end{lemma}

Now let $(M,\omega)$ be a $2n$-dimensional compact symplectic
manifold.  Applying the above construction to each tangent space we
obtain a symplectic Hodge star operator
$*\colon\Omega^k(M)\to\Omega^{2n-k}(M)$ satisfying
$$
\alpha\wedge*\beta=\omega(\alpha,\beta)\frac{\omega^n}{n!}
$$
for $\alpha$, $\beta\in\Omega^k(M)$.  The
symplectic boundary operator $\delta\colon\Omega^k\to\Omega^{k-1}$ is
then defined by 
$$\delta\alpha=(-1)^{k+1}{*}d{*}\alpha$$
for $\alpha\in\Omega^k(M)$.  As in Riemannian Hodge theory we call
$\alpha$ \emph{coclosed} if $\delta\alpha=0$, \emph{coexact} if
$\alpha=\delta\beta$ for some $\beta$, and \emph{harmonic} if it is
closed and coclosed.  However, in a striking departure from the
Riemannian case, the symplectic Laplacian $d\delta+\delta d$ vanishes
identically, i.e.\ $d$ and $\delta $ anti-commute.

\begin{example}\label{example;coexact}
Every function is coclosed.  A function $f$ is coexact if and only if
the $2n$-form $*f=f\,{*}1=f\omega^n/n!$ is exact, which is the case if
and only if $\int_Mf\omega^n/n!=0$, i.e.\ the mean of $f$ relative to
the Liouville measure is zero.
\end{example}

Let $\pi=\omega^\sharp$ be the Poisson two-vector defined by the
symplectic form and let $v_f=(df)^\sharp$ denote the Hamiltonian
vector field of a function $f$.  The first part of the following
result, which shows that $\delta$ can be defined on any Poisson
manifold, is Koszul's original definition of $\delta$ (see
\cite{koszul;crochet}) and the second part is a Leibniz rule for
$\delta$.

\begin{proposition}\label{proposition;delta}
\begin{enumerate}
\item\label{item;brylinski}
$\delta=[\iota(\pi),d]$.
\item\label{item;leibniz}
$\delta(f\alpha)=f\delta\alpha-\iota(v_f)\alpha$ for all functions $f$
and forms $\alpha$. 
\item\label{item;coexact}
If $\alpha$ is closed, then $\delta\alpha$ is exact.  If $\alpha$ is
coclosed, then $d\alpha$ and $\iota(v_f)\alpha$ are coexact.
\end{enumerate}
\end{proposition}

\begin{proof}
\eqref{item;brylinski} is due to
\cite[Theorem~2.2.1]{brylinski;differential-poisson}.  Let $\eps(f)$
denote multiplication by $f$.  Using \eqref{item;brylinski} and the
ordinary Leibniz rule $[d,\eps(f)]=\eps(df)$ we get
\begin{equation}\label{equation;delta-eps}
[\delta,\eps(f)]=[[\iota(\pi),d],\eps(f)]
=[\iota(\pi),[d,\eps(f)]]-[d,[\iota(\pi),\eps(f)]]=[\iota(\pi),\eps(df)].
\end{equation}
The identity  $df=\iota(v_f)\omega$ yields
$\eps(df)=[\iota(v_f),\eps(\omega)]$, so
\begin{multline*}
[\iota(\pi),\eps(df)]=[\iota(\pi),[\iota(v_f),\eps(\omega)]]
=[\iota(v_f),[\iota(\pi),\eps(\omega)]]
+[[\iota(\pi),\iota(v_f)],\eps(\omega)]\\
=[\iota(v_f),A]=-\iota(v_f),
\end{multline*}
where we used $[\iota(\pi),\eps(\omega)]=A$, the scalar operator
defined by $A\alpha=(n-k)\alpha$ for $\alpha$ of degree $k$.  Together
with \eqref{equation;delta-eps} this shows
$[\delta,\eps(f)]=-\iota(v_f)$, which is equivalent to
\eqref{item;leibniz}.  The first statement in \eqref{item;coexact}
follows immediately from \eqref{item;brylinski}.  (See
\cite[Proposition 1.3.2]{brylinski;differential-poisson}.)  Applying
Hodge star we get that $d\alpha$ is coexact for $\alpha$ coclosed.
Finally, if $\delta\alpha=0$, then $\iota(v_f)\alpha=-\delta(f\alpha)$
by \eqref{item;leibniz}.
\end{proof}

The following result of Mathieu can be thought of as an analogue of
the Hodge theorem.  (See also \cite{yan;hodge-structure-symplectic}
for an alternative proof.)  Let $H(M)=H(\Omega(M))$ denote the de Rham
cohomology of $M$.

\begin{theorem}[\cite{mathieu;harmonic-symplectic}]
\label{theorem;mathieu}
The following conditions are equivalent.
\begin{enumerate}
\item
$M$ satisfies the strong Lefschetz condition, i.e.\ the map
$$
H^{n-k}(M)\longrightarrow H^{n+k}(M),\qquad
c\longmapsto[\omega]^k\wedge c,
$$
is an isomorphism for each $0\leq k\leq n$.
\item
Every cohomology class in $H(M)$ has a harmonic representative.
\end{enumerate}
\end{theorem}

Let $\Omega_\delta(M)$ be the kernel of $\delta$.  Since $d$ and
$\delta$ anti-commute, $\Omega_\delta(M)$ forms a subcomplex of
$\Omega(M)$, the cohomology of which we denote by $H_\delta(M)$.  In
this language Mathieu's result can be restated as follows.

\begin{theorem}\label{theorem;mathieu-surjective}
The strong Lefschetz property holds for $(M,\omega)$ if and only if
the map $H_\delta(M)\to H(M)$ induced by the inclusion
$\Omega_\delta(M)\inj\Omega(M)$ is surjective.
\end{theorem}

This restatement is useful for the following reason.  The identity
$$\ker d\cap\im\delta=\im d\delta=\ker\delta\cap\im d$$
is known as the \emph{$d\delta$-lemma}.  In plain language: if a
$k$-form $\alpha$ is harmonic and either exact or coexact, then
$\alpha=d\delta\beta$ for some $k$-form $\beta$.  (In particular
$\alpha$ is both exact \emph{and} coexact.)  If the $d\delta$-lemma
holds for $M$, then for every closed form $\gamma$ the equation
$\delta\gamma=d\delta\zeta$ is soluble for $\zeta$ and so the
cohomology class of $\gamma$ has a harmonic representative.  Therefore
$M$ satisfies strong Lefschetz by Theorem
\ref{theorem;mathieu-surjective}.  Intriguingly, the converse is also
true.

\begin{theorem}%
[\cite{merkulov;formality-canonical-symplectic,guillemin;symplectic-hodge}]
\label{theorem;ddelta}
Strong Lefschetz holds if and only if the $d\delta$-lemma holds.
\end{theorem}

The following result is an easy consequence.  Here
$H(\Omega(M),\delta)$ denotes the homology of $\Omega(M)$ with respect
to $\delta$.

\begin{theorem}%
[formality,
\cite{merkulov;formality-canonical-symplectic,guillemin;symplectic-hodge}]
\label{theorem;quasi}
Assume strong Lefschetz.  Then the $d$-chain maps in the diagram
$$\Omega(M)\longleftarrow\Omega_\delta(M)\longrightarrow
H(\Omega(M),\delta)$$
are quasi-isomorphisms, i.e.\ they induce isomorphisms in cohomology.
Therefore the de Rham complex of $M$ is formal, i.e.\ quasi-isomorphic
to a complex with zero differential.
\end{theorem}

This immediately implies the following sharpened version of Theorem
\ref{theorem;mathieu-surjective}.

\begin{theorem}[\cite{guillemin;symplectic-hodge}]
\label{theorem;mathieu-bijective}
Strong Lefschetz holds if and only if the map $H_\delta(M)\to H(M)$ is
bijective.
\end{theorem}

%%%%%%%%%%%%%%%%%%%%%%%%%%%%%%%%%%%%%%%%%%%%%%%%%%%%%%%%%%%%%%%%%%%%%%%%
\section{Equivariant formality and the $d_G\delta$-lemma}
\label{section;dgdelta}
%%%%%%%%%%%%%%%%%%%%%%%%%%%%%%%%%%%%%%%%%%%%%%%%%%%%%%%%%%%%%%%%%%%%%%%%

In this section we establish equivariant versions of Theorems
\ref{theorem;ddelta}--\ref{theorem;quasi}.  Let $G$ be a connected
compact Lie group acting in a Hamiltonian fashion on a
$2n$-dimensional compact symplectic manifold $(M,\omega)$.  Let
$\phi\colon\g\to C^\infty(M)$ be an equivariant moment map.
Throughout this section we assume $(M,\omega)$ to have the strong
Lefschetz property.  (It is well-known that the Lefschetz and
Hamiltonian conditions are not independent.  For instance, under the
Lefschetz assumption every symplectic $G$-action whose generating
vector fields have fixed points is Hamiltonian.)

We start with a quick review of equivariant de Rham theory.  See e.g.\
\cite{guillemin-sternberg;supersymmetry-equivariant} for more details.
Let $\Omega_G(M)=(S\g^*\otimes\Omega(M))^G$ be the Cartan double
complex of the $G$-manifold $M$.  For brevity let us write
$\Omega=\Omega(M)$ and $\Omega_G=\Omega_G(M)$.  Elements of $\Omega_G$
can be regarded as equivariant polynomial maps from $\g$ to $\Omega$
and are called \emph{equivariant differential forms} on $M$.  The
bigrading on the Cartan complex is defined by
$\Omega_G^{ij}=(S^i\g^*\otimes\Omega^{j-i})^G$.  It is equipped with
the vertical differential $1\otimes d$, which we will abbreviate to
$d$, and the horizontal differential $\partial$, which is defined by
$\partial\alpha(\xi)=-\iota(\xi)\alpha(\xi)$.  Here $\iota(\xi)$
denotes inner product with the vector field on $M$ induced by
$\xi\in\g$.  We can regard $\Omega_G$ as a single complex with grading
$\Omega_G^k=\bigoplus_{i+j=k}\Omega_G^{ij}$ and total differential
$d_G=d+\partial$, which is called the \emph{equivariant exterior
derivative}.  If $d_G\alpha=0$, resp.\ $\alpha=d_G\beta$, we call
$\alpha$ \emph{equivariantly closed}, resp.\ \emph{equivariantly
exact}.  The total cohomology $\ker d_G/\im d_G$ is the
\emph{equivariant de Rham cohomology} $H_G(M)$.  Observe that
$\Omega_G^{0j}=(\Omega^G)^j$, the space of \emph{in}variant $j$-forms
on $M$.  Thus the zeroth column of the Cartan complex $\Omega_G$ is
the invariant de Rham complex $\Omega^G$, which is a deformation
retract of the ordinary de Rham complex $\Omega$ because $G$ is
connected.  Hence $H(\Omega^G)=H(\Omega)=H(M)$.  The Cartesian
projection $\bar{p}\colon\Omega_G\to\Omega^G$, defined by
$\bar{p}(\alpha)=\alpha(0)$, is a chain map with respect to the
equivariant exterior derivative $d_G$ on $\Omega_G$ and the ordinary
exterior derivative $d$ on $\Omega^G$.  The Cartan complex is a de
Rham model for the homotopy quotient $M_G$, which is the total space
of the bundle with fibre $M$ associated to the universal $G$-bundle
over the classifying space $B_G$.  The projection $\bar{p}$ models the
inclusion of the fibre $M\inj M_G$.

However, symplectic Hodge theory gives us a fourth differential
$1\otimes\delta$ of bidegree $(0,-1)$, which we will write as
$\delta$.

\begin{lemma}\label{lemma;partial-delta}
$\partial\delta=-\delta\partial$ and $d_G\delta=-\delta d_G$.
\end{lemma}

\begin{proof}
Let $\pi$ be the Poisson two-vector defined by the symplectic form.
Using Proposition \ref{proposition;delta}\eqref{item;brylinski} we
find for all $\alpha\in\Omega_G$
\begin{align*}
\partial\delta\alpha(\xi)
&=\bigl(-\iota(\xi)\iota(\pi)d+\iota(\xi)d\iota(\pi)\bigr)\alpha(\xi)\\
&=\bigl(\iota(\pi)d\iota(\xi)-d\iota(\pi)\iota(\xi)
-\iota(\pi)\ca{L}(\xi)+\ca{L}(\xi)\iota(\pi)\bigr)\alpha(\xi),\\
\delta\partial\alpha(\xi)
&=\bigl(-\iota(\pi)d\iota(\xi)+d\iota(\pi)\iota(\xi)\bigr)\alpha(\xi),
\end{align*}
so that $(\partial\delta+\delta\partial)\alpha(\xi)
=\bigl(\ca{L}(\xi)\iota(\pi)-\iota(\pi)\ca{L}(\xi)\bigr)\alpha(\xi)$.
Now $\alpha\bigl(\Ad g^{-1}(\xi)\bigr)=g^*\alpha(\xi)$ for all $g\in
G$ and $\xi\in\g$ and hence
$\alpha([\eta,\xi])=\ca{L}(\eta)\alpha(\xi)$ for all $\xi$,
$\eta\in\g$.  Likewise, since $\pi$ is invariant,
$\iota(\pi)\alpha\bigl(\Ad g^{-1}(\xi)\bigr)
=g^*\iota(\pi)\alpha(\xi)$ and so
$\iota(\pi)\alpha([\eta,\xi])=\ca{L}(\eta)\iota(\pi)\alpha(\xi)$.
Hence $\iota(\pi)\ca{L}(\xi)\alpha(\xi)
=\ca{L}(\xi)\iota(\pi)\alpha(\xi)=0$ and
$(\partial\delta+\delta\partial)\alpha(\xi)=0$, i.e.\
$\partial\delta=-\delta\partial$.  Since $d\delta=-\delta d$, this
implies $d_G\delta=(d+\partial)\delta=-\delta(d+\partial)=-\delta
d_G$.
\end{proof}

This implies that $\Omega_{G,\delta}=\ker\delta\cap\Omega_G$ is a
double subcomplex of $\Omega_G$ and that the homology
$H(\Omega_G,\delta)$ of $\Omega_G$ with respect to $\delta$ is a
double complex with differentials induced by $d$ and $\partial$.  Thus
we have a diagram of morphisms of double complexes
\begin{equation}\label{equation;double}
\Omega_G\longleftarrow\Omega_{G,\delta}\longrightarrow
H(\Omega_G,\delta).
\end{equation}
Since $\delta$ does not act on the polynomial part, these morphisms
are linear over the invariant polynomials $(S\g^*)^G$.  Let us first
examine the complex $H(\Omega_G,\delta)$.  We need a preliminary
result about the action of $\iota(\xi)$ on invariant forms.  It
follows from the identity $\ca{L}(\xi)=d\iota(\xi)+\iota(\xi)d$ that
$\iota(\xi)\colon\Omega^G\to\Omega^G$ is a chain map (of degree $1$)
with respect to $d$.  Likewise, the identity
$\partial\delta+\delta\partial=0$ of Lemma \ref{lemma;partial-delta}
applied to the zeroth column of $\Omega_G$ shows that $\iota(\xi)$ is
a chain map with respect to $\delta$.

\begin{lemma}\label{lemma;ginzburg}
Let $\xi\in\g$ and $\alpha\in\Omega^G$.  If $\alpha$ is closed, then
$\iota(\xi)\alpha$ is exact.  If $\alpha$ is coclosed, then
$\iota(\xi)\alpha$ is coexact.  Equivalently, for every $\xi\in\g$ the
maps
$$
\iota_d\colon H(\Omega^G,d)\to H(\Omega^G,d),\qquad\iota_\delta\colon
H(\Omega^G,\delta)\to H(\Omega^G,\delta)
$$
induced by $\iota(\xi)\colon\Omega^G\to\Omega^G$ are zero.
\end{lemma}

\begin{proof}
The fact that $\iota_\delta=0$ (which is true without the Lefschetz
condition) follows from Proposition
\ref{proposition;delta}\eqref{item;coexact} and the assumption that
the $G$-action is Hamiltonian.  The assertion $\iota_d=0$ is due to
Ginzburg \cite{ginzburg;equivariant-cohomology-kahler}.  For the
reader's convenience we recall the argument.  Let
$\alpha\in\Omega^{n-k}$ be invariant and closed.  Recall that the
class $[\alpha]\in H^{n-k}(M)$ is called \emph{primitive} if
$[\omega^{k+1}\wedge\alpha]=0$, and that strong Lefschetz implies the
Lefschetz decomposition $H^i(M)=\bigoplus_{j\ge0}[\omega^j]P^{i-2j}$,
where $P^i\subset H^i(M)$ is the subspace of primitive elements.
First consider the case that $[\alpha]$ is primitive.  Suppose
$\iota_d[\alpha]\in H^{n-k+1}(M)$ was nonzero.  The identity
$\iota(\xi)(\omega\wedge\alpha)
=d\phi(\xi)\wedge\alpha+\omega\wedge\iota(\xi)\alpha$ implies that
$\iota_d$ commutes with multiplication by the symplectic class
$[\omega]$.  Thus by strong Lefschetz
$$
\iota_d\bigl[\omega^{k+1}\wedge\alpha\bigr]
=[\omega]^{k+1}\wedge\iota_d[\alpha]\ne0,
$$
which contradicts $[\omega^{k+1}\wedge\alpha]=0$.  Therefore
$\iota_d[\alpha]=0$ for $[\alpha]$ primitive.  It now follows from the
primitive decomposition that $\iota_d[\alpha]=0$ for all $[\alpha]$.
\end{proof}

This enables us to show that $H(\Omega_G,\delta)$ is a trivial double
complex.

\begin{lemma}\label{lemma;zero}
Both differentials on $H(\Omega_G,\delta)$ are zero.
\end{lemma}

\begin{proof}
First note that the (ordinary) $d\delta$-lemma, Theorem
\ref{theorem;ddelta}, holds for equivariant forms as well as for
ordinary forms.  The reason is that $d$ and $\delta$ act on $\Omega_G$
as $1\otimes d$ and $1\otimes\delta$ and that both operators are
$G$-equivariant.  Now suppose $\alpha\in\Omega_G$ satisfies
$\delta\alpha=0$.  Then $d\alpha=d\delta\beta$ for some
$\beta\in\Omega_G$ by the $d\delta$-lemma.  Hence $d\alpha=-\delta
d\beta$, so the differential on $H(\Omega_G,\delta)$ induced by $d$ is
zero.  To prove that the other differential is zero we must be more
careful about picking a representative of an element of
$H(\Omega_G,\delta)$.  Observe that
\begin{equation}\label{equation;homology}
H(\Omega_G,\delta)=(S\g^*\otimes H(\Omega,\delta))^G=(S\g^*)^G\otimes
H(M),
\end{equation}
where we used the isomorphism $H(\Omega,\delta)\cong H(M)$ of Theorem
\ref{theorem;quasi} and the connectedness of $G$ to ensure that $G$
acts trivially on $H(M)$ and $H(\Omega,\delta)$.  (It is easy to see
that the inclusion $\Omega^G\inj\Omega$ is a deformation retraction
for $\delta$ as well as for $d$.)  Choose a basis $\{\xi_q\}$ of $\g$,
let $\{x_q\}$ be the dual basis of $\g^*$ and let $\{f_p\}$ be a basis
of the vector space $(S\g^*)^G$ of invariant polynomials.  It follows
from \eqref{equation;homology} that an element of $H(\Omega_G,\delta)$
can be represented by an $\alpha\in\Omega_G$ with $\delta\alpha=0$ of
the form $\alpha=\sum_pf_p\otimes\alpha_p$ for unique
$\alpha_p\in\Omega^G$.  Then we have $\delta\alpha_p=0$ for all $p$.
Using Lemma \ref{lemma;ginzburg} we find $\beta_{pq}\in\Omega^G$ such
that $\iota(\xi_q)\alpha_p=d\beta_{pq}$.  Thus
$$
\partial\alpha=\sum_{pq}x_qf_p\otimes\iota(\xi_q)\alpha_p=d\beta
$$
where $\beta=\sum_{pq}x_qf_p\otimes\beta_{pq}\in S\g^*\otimes\Omega$.
After averaging over $G$ we get $\partial\alpha=d\beta$ with
$\beta\in\Omega_G$.  Furthermore
$\delta\partial\alpha=-\partial\delta\alpha=0$, so $\partial\alpha$ is
exact and coclosed.  Applying once again the $d\delta$-lemma we see
that $\partial\alpha$ is coexact, i.e.\ the differential on
$H(\Omega_G,\delta)$ induced by $\partial$ is zero.
\end{proof}

Let $E$ be the spectral sequence of $\Omega_G$ relative to the
filtration associated to the horizontal grading and $E_\delta$ that of
$\Omega_{G,\delta}$.  The first terms are
\begin{equation}\label{equation;spectral}
\begin{split}
E_1&=\ker d/\im d=(S\g^*\otimes H(M))^G=(S\g^*)^G\otimes H(M),\\
(E_\delta)_1&=(\ker d\cap\ker\delta)/(\im
d\cap\ker\delta)=(S\g^*\otimes H_\delta(M))^G=(S\g^*)^G\otimes H(M).
\end{split}
\end{equation}
Here we used the isomorphism $H_\delta(M)\cong H(M)$ of Theorem
\ref{theorem;quasi} and the connectedness of $G$.  By Lemma
\ref{lemma;zero} $H(\Omega_G,\delta)$ is a trivial complex, so its
spectral sequence is constant with trivial differentials at each
stage.  The two morphisms \eqref{equation;double} induce morphisms of
spectral sequences
$$
E\longleftarrow E_\delta\longrightarrow H(\Omega_G,\delta).
$$
It follows from \eqref{equation;spectral} and
\eqref{equation;homology} that these morphisms are isomorphisms at the
first stage.  Hence they are isomorphisms at \emph{every} stage.  In
particular the three sequences converge to the same limit, so the
morphisms \eqref{equation;double} induce isomorphisms on total
cohomology.  In fact, since the spectral sequence for
$H(\Omega_G,\delta)$ is constant, so are the spectral sequences $E$
and $E_\delta$.  This proves the following result, where
$H_{G,\delta}(M)$ denotes the total cohomology of $\Omega_{G,\delta}$.

\begin{theorem}[equivariant formality]\label{theorem;quasi-equivariant}
The spectral sequences $E$ and $E_\delta$ degenerate at the first
term.  The morphisms \eqref{equation;double} induce isomorphisms of
$(S\g^*)^G$-modules
$$
H_G(M)\overset\cong{\longleftarrow}H_{G,\delta}(M)
\overset\cong{\longrightarrow}(S\g^*)^G\otimes H(M).
$$
Hence the equivariant de Rham complex of $M$ is formal, i.e.\
quasi-isomorphic to a complex with zero differential.
\end{theorem}

The spectral sequence $E$ is isomorphic to the Leray spectral sequence
of the fibre bundle $M_G\to B_G$ and its degeneracy means that the map
induced by $\bar{p}$,
$$
p\colon H_G(M)\to H(M),
$$
is surjective, in other words that each cohomology class on the fibre
$M$ can be extended to an equivariant cohomology class.  We assert
that there is a natural choice of such an extension.  Let
$$
s\colon H(M)\to H_G(M)
$$
be the composition of the map $H(M)\inj(S\g^*)^G\otimes H(M)$ which
sends a cohomology class $c$ to $1\otimes c$ and the isomorphism
$(S\g^*)^G\otimes H(M)\to H_G(M)$ given by Theorem
\ref{theorem;quasi-equivariant}.

\begin{corollary}\label{corollary;extension}
$s$ is a section of $p$.  Thus every cohomology class has a canonical
equivariant extension.
\end{corollary}

\begin{proof}
Via the isomorphisms $H(M)\cong H(\Omega,\delta)$ and $H_G(M)\cong
H(\Omega_G,\delta)$ given by Theorems \ref{theorem;quasi} and
\ref{theorem;quasi-equivariant} the map $p$ corresponds to the natural
projection $(S\g^*)^G\otimes H(\Omega,\delta)\to H(\Omega,\delta)$ and
$s$ corresponds to the inclusion $H(\Omega,\delta)\inj(S\g^*)^G\otimes
H(\Omega,\delta)$ defined by $c\mapsto1\otimes c$.  Therefore $p\circ
s=1$.
\end{proof}

\begin{corollary}\label{corollary;spectral}
If $\alpha\in\Omega_G$ is closed, then $\partial\alpha$ is exact.
More generally, let $j\ge1$ and suppose for $0\le i<j$ equivariant
forms $\zeta_i\in\Omega_G$ are given such that $d\zeta_0=0$ and
$\partial\zeta_{i-1}+d\zeta_i=0$ for $1\le i<j$.  Then there exists
$\zeta_j\in\Omega_G$ such that $\partial\zeta_{j-1}+d\zeta_j=0$.
\end{corollary}

\begin{proof}
This is equivalent to the degeneracy of $E$.  (Cf.\- the description of
the differentials on $E$ in
\cite[p.~164]{bott-tu;differential-forms}.)
\end{proof}

These results improve on
\cite[Proposition~5.8]{kirwan;cohomology-quotients-symplectic} and
\cite[Theorem 3.3]{ginzburg;equivariant-cohomology-kahler}, which
state that $E$ degenerates and that $H_G(M)\cong S\g^*\otimes H(M)$
non-canonically.  However, the cited results hold for arbitrary
compact Hamiltonian $G$-manifolds, whereas we have imposed the
Lefschetz condition.  For completeness we record how in principle one
can compute the section $s$.  Let $[\alpha]\in H^k(M)$.  Let us
assume, as we may, that the representative $\alpha$ is invariant and
harmonic.  Then $d\alpha=0$ implies that $\partial\alpha$ is exact by
the first statement in Corollary \ref{corollary;spectral}, and
$\delta\alpha=0$ implies that $\partial\alpha$ is coclosed.  The
$d\delta$-lemma gives $\beta_1\in\Omega_G^{1k}$ such that
$$
\partial\alpha=-d\delta\beta_1.
$$
This equation shows that $\partial d\delta\beta_1=0$, so
$d\partial\delta\beta_1=-\partial d\delta\beta_1=0$.  Also
$\partial\delta\beta_1=-\delta\partial\beta_1$ is coexact, so another
application of the $d\delta$-lemma yields
$\partial\delta\beta_1=-d\delta\beta_2$ for some
$\beta_2\in\Omega_G^{2,k-1}$.  Recursively solving the equation
$\partial\delta\beta_i=-d\delta\beta_{i+1}$ for $\beta_{i+1}$ in this
manner we find a finite chain of $\beta_i\in\Omega_G^{i,k+1-i}$.  Put
$$\alpha_G=\alpha+\sum_i\delta\beta_i\in\Omega_G^k.$$
Then $d_G\alpha_G=\delta\alpha_G=0$ and $\bar{p}\alpha_G=\alpha$ by
construction.  Let $\{\alpha\}$, resp.\ $\{\alpha_G\}$, be the
homology class of $\alpha$ in $H(\Omega,\delta)$, resp.\ of $\alpha_G$
in $H(\Omega_G,\delta)=(S\g^*)^G\otimes H(\Omega,\delta)$.  Since
$\alpha_G-\alpha=\sum_i\delta\beta_i$ is coexact, the element
$1\otimes\{\alpha\}\in(S\g^*)^G\otimes H(\Omega,\delta)$ is equal to
$\{\alpha_G\}$.  Therefore $s([\alpha])=[\alpha_G]$.  (The point of
Corollary \ref{corollary;extension} is that the equivariant cohomology
class of $\alpha_G$ is independent of the choice of the $\beta_i$.)

\begin{example}\label{example;extension-omega}
The equivariant symplectic form $\omega+\phi\in\Omega_G^2$ is an
equivariantly closed extension of the symplectic form, but $\phi$ is
not necessarily coexact, so we may not have $s[\omega]=[\omega+\phi]$.
Put $\chi_0(\xi)=\int_M\phi(\xi)\omega^n/n!$.  Then
$\chi_0\colon\g\to\R$ is a character of $\g$.  According to Example
\ref{example;coexact} the shifted moment map $\phi_0=\phi-\chi_0$ is
coexact, so $\omega_G=\omega+\phi_0$ is the equivariant extension for
which $s[\omega]=[\omega_G]$.
\end{example}

As an aside we point out that $s$ (or any other section of $p$) is
very far from being a ring homomorphism.  Indeed, if $s$ was
multiplicative we would have
$0=s\bigl(\omega^{n+1}\bigr)=[\omega_G]^{n+1}$, i.e.\
$(\omega_G)^{n+1}=\phi\omega^n+\cdots+\phi^{n+1}$ would be
equivariantly exact.  In particular $\partial\lambda=\phi^{n+1}$ for
some $\lambda\in\Omega_G^{n,n+1}$, i.e.\
$-\iota(\xi)\lambda(\xi)=\phi(\xi)^{n+1}$ for all $\xi\in\g$.  The
extrema of $\phi(\xi)$ are critical points for the vector field
induced by $\xi$, so this equality implies that $\phi(\xi)$ vanishes
identically.  Thus $s$ is multiplicative only in the uninteresting
case of a trivial action.

\begin{example}
As the preceding remark indicates, extending the powers of the
symplectic class is not entirely straightforward.  As an example we
find $s\bigl([\omega]^2\bigr)$.  The moment map $\phi_0$ of Example
\ref{example;extension-omega} is coexact, so $\phi_0=\delta\phi_1$ for
some $\phi_1\in\Omega_G^{12}$.  Then
$\frac12\phi_0^2+\partial\phi_1\in\Omega_G^{22}$, so
$$
\chi_1(\xi)=\int_M\Bigl(\frac12\phi_0(\xi)^2+\partial\phi_1(\xi)\Bigr)
\frac{\omega^n}{n!}
$$
defines an invariant quadratic polynomial on $\g$.  Moreover, by
Example \ref{example;coexact}
$\frac12\phi_0^2+\partial\phi_1-\chi_1=\delta\phi_2$ for some
$\phi_2\in\Omega_G^{23}$.  A computation shows that
$\omega^2-2\delta(\omega\wedge\phi_1)-\delta\phi_2$ is equivariantly
closed.  Therefore
$$
s\bigl([\omega]^2\bigr)
=\bigl[\omega^2-2\delta(\omega\wedge\phi_1)-\delta\phi_2\bigr].
$$
\end{example}

A variation on the construction of equivariant extensions gives us the
$d_G\delta$-lemma.  Call an equivariant form $\alpha$
\emph{equivariantly harmonic} if it is equivariantly closed
($d_G\alpha=0$) and coclosed ($\delta\alpha=0$).

\begin{theorem}[$d_G\delta$-lemma]\label{theorem;dgdelta}
Suppose $\alpha\in\Omega_G$ is equivariantly harmonic and either
equivariantly exact or coexact.  Then there exists $\beta\in\Omega_G$
such that $\alpha=d_G\delta\beta$.
\end{theorem}

\begin{proof}
Let $\alpha\in\Omega_G^k$ and assume $\alpha=d_G\gamma$ and
$\delta\alpha=0$.  Decomposing $\alpha=\sum_i\alpha_i$ and
$\gamma=\sum_i\gamma_i$ into homogeneous parts
$\alpha_i\in\Omega_G^{i,k-i}$ and $\gamma_i\in\Omega_G^{i,k-i-1}$ we
have
\begin{align}
\delta\alpha_i&=0\tag{$*_i$}\label{equation;*i}\\
\alpha_i&=\partial\gamma_{i-1}+d\gamma_i
\tag{$\dagger_i$}\label{equation;daggeri}
\end{align}
for $i\ge0$, where we put $\gamma_{-1}=0$.  We need to solve
$\alpha=d_G\delta\beta$ for $\beta\in\Omega^k_G$, i.e.\ we must find
$\beta_i\in\Omega_G^{i,k-i}$ such that
\begin{equation}\label{equation;ddaggeri}
\alpha_i=\partial\delta\beta_{i-1}+d\delta\beta_i\tag{$\ddagger_i$}
\end{equation}
for $i\ge0$, where again we put $\beta_{-1}=0$.  The zeroth equation
($\ddagger_0$) amounts to $\alpha_0=d\delta\beta_0$, which is soluble
by the $d\delta$-lemma because the integrability conditions ($*_0$)
and ($\dagger_0$) say that $\alpha_0$ is coclosed and exact.  Now let
$j\ge1$ and assume that \eqref{equation;ddaggeri} has been solved for
$i<j$.  Write $\zeta_i=\gamma_i-\delta\beta_i$ for $-1\le i<j$.  Then
$\zeta_{-1}=0$.  Also
$$
d\zeta_i=d(\gamma_i-\delta\beta_i)
=-\partial\gamma_{i-1}+\partial\delta\beta_{i-1}=-\partial\zeta_{i-1}
$$
for $0\le i<j$ by \eqref{equation;daggeri} and
\eqref{equation;ddaggeri}.  Hence by Corollary
\ref{corollary;spectral} there exists $\zeta_j$ such that
$\partial\zeta_{j-1}+d\zeta_j=0$.  Now rewrite ($\ddagger_j$) as
$$
d\delta\beta_j=\alpha_j-\partial\delta\beta_{j-1}
$$
and observe that by ($*_j$) the right-hand side of this equation is
coclosed.  Using ($\dagger_j$) we rewrite the right-hand side as
$$
\alpha_j-\partial\delta\beta_{j-1}
=\partial\gamma_{j-1}+d\gamma_j-\partial\delta\beta_{j-1}
=\partial\zeta_{j-1}+d\gamma_j=d(-\zeta_j+\gamma_j),
$$
which is exact.  So ($\ddagger_j$) is soluble by the $d\delta$-lemma.
By induction \eqref{equation;ddaggeri} is soluble for all $i\ge0$, so
we find $\beta$ with $\alpha=d_G\delta\beta$.  The second half of the
$d_G\delta$-lemma follows from the first: assume $d_G\alpha=0$ and
$\alpha$ coexact.  Then the homology class of $\alpha$ in
$H(\Omega_G,\delta)$ is zero, so by Theorem
\ref{theorem;quasi-equivariant} the cohomology class of $\alpha$ in
$H_{G,\delta}(M)$ is zero, i.e.\ $\alpha$ is equivariantly exact.
Hence $\alpha=d_G\delta\beta$ for some $\beta$.
\end{proof}

Let $\Omega_{G,\mathrm{per}}^{ij}=\Omega_G^{j-i}$ be the
\emph{periodic} equivariant de Rham complex of $M$, equipped with the
horizontal differential $\delta$ and the vertical differential $d_G$.
The $d_G\delta$-lemma has as an immediate consequence the following
equivariant analogue of \cite[Theorem
2.3.1]{brylinski;differential-poisson}.

\begin{corollary}\label{corollary;periodic}
The spectral sequences of $\Omega_{G,\mathrm{per}}$ with respect to
the horizontal and vertical filtrations both degenerate at the first
term.
\end{corollary}

See \cite[Section 2.3]{brylinski;differential-poisson} for the
definition of degeneration for a double complex that is zero below the
diagonal $i=j$.  (But note Brylinski's opposite conventions for the
periodic complex.)  As in \emph{loc.\ cit}., the degeneracy with
respect to the vertical filtration is probably true even if $M$ does
not satisfy the Lefschetz condition.

%%%%%%%%%%%%%%%%%%%%%%%%%%%%%%%%%%%%%%%%%%%%%%%%%%%%%%%%%%%%%%%%%%%%%%%%
%%%%%%%%%%%%%%%%%%%%%%%%%%%%%%%%%%%%%%%%%%%%%%%%%%%%%%%%%%%%%%%%%%%%%%%%

\providecommand{\bysame}{\leavevmode\hbox to3em{\hrulefill}\thinspace}
\providecommand{\MR}{\relax\ifhmode\unskip\space\fi MR }
\providecommand{\MRhref}[2]{%
  \href{http://www.ams.org/mathscinet-getitem?mr=#1}{#2}
}
\providecommand{\href}[2]{#2}

%%%%%%%%%%%%%%%%%%%%%%%%%%%%%%%%%%%%%%%%%%%%%%%%%%%%%%%%%%%%%%%%%%%%%%%%
%%%%%%%%%%%%%%%%%%%%%%%%%%%%%%%%%%%%%%%%%%%%%%%%%%%%%%%%%%%%%%%%%%%%%%%%

\end{document}